# A New Problem in String Searching

George Havas* and Jin Xian Lian**

Key Centre for Software Technology
Department of Computer Science
The University of Queensland
Queensland 4072, Australia

**Abstract.** We describe a substring search problem that arises in group presentation simplification processes. We suggest a two-level searching model: skip and match levels. We give two timestamp algorithms which skip searching parts of the text where there are no matches at all and prove their correctness. At the match level, we consider Harrison signature, Karp-Rabin fingerprint, Bloom filter and automata based matching algorithms and present experimental performance figures.

## 1 Introduction

A fundamental technique used in computer science is to search for a specific substring in a large body of text. Text-processing systems must allow their users to search for given character strings within a body of text. Database systems must be capable of searching for records with stated values in specified fields. Substring searching plays an important role in group presentation simplification processes too. Most of execution time of these processes is used in the substring search part.

We describe the substring search problem that arises in presentation simplification processes, also known as Tietze processes, and give two timestamp algorithms which skip searching parts of the text where there are no matches at all. First we give the background of Tietze processes then we formalize the substring searching problem in this context and define the terminology used in this paper. We analyze the problem and give a two-level searching model. Based on this analysis, we study both the skip and match levels. We present two timestamp algorithms at the skip level and two minimal-cover theorems for those algorithms. At the match level, we consider algorithms and data structures based on Harrison signatures, Karp-Rabin fingerprints, Bloom filters and automata. We indicate the practical performance of these in this context.

Finitely presented groups have been much studied. All requisite mathematical background is provided in [14, Chapter 1]. An overview of algorithms for such groups is included in [5], and a comprehensive book on computation with finitely presented groups [19] has recently appeared.

---

* email: havas@cs.uq.oz.au; partially supported by the Australian Research Council.
** email: lili@cs.uq.oz.au; partially supported by an Australian Overseas Postgraduate Research Scholarship.

A finitely presented group may be given by a presentation $G = \langle g_1, \ldots, g_d \mid R_1, \ldots, R_n \rangle$ where the $g_i$ are generators and the $R_j$ are relators. Generally speaking, presentations are good if they are short: few generators; and few relators of reasonable length. This makes them relatively intelligible to humans and also often makes them better suited for computer calculations. We are interested in the situation where we have what we regard as a bad presentation for a group and we wish to find a good presentation. This kind of situation may arise in a number of ways.

A theorem of Tietze proves that, given two presentations of a group, there exists a sequence of simple transformations which demonstrates that the presentations are of the same group. However there is no general algorithm for finding such a sequence, a consequence of unsolvability results in group theory. Various Tietze transformation procedures, which input a "bad" presentation and output a "good" presentation for a group, have been described ([9, 17, 11]). Newer procedures, written in the higher level GAP [18] language (with special kernel support), have also been developed by Volkmar Felsch and Martin Schönert in Aachen. Three main principles used by Tietze transformation methods to simplify presentations are: short eliminations; long eliminations; and substring replacements.

In each short elimination phase, all relators of length 1 and non-involutory relators of length 2 are used to eliminate generators and their associated relators. In each long elimination phase, redundant generators (generators which occur only once in some relator) and their associated relators are eliminated using relators with length greater than 2.

In each substring replacement phase, relators are shortened by replacing long substrings with shorter equivalent strings. First substring searching is performed. A relator $R_i$ is chosen and other relators are searched for a matching substring $v$ in a rotation $uv$ of $R_i$ or its inverse and in a rotation $wv$ of $R_j$, with the length of $v$ greater than the length of $u$. Then, when such a match is found, the relator $R_j$ is replaced by the shorter relator $wu^{-1}$. One substring replacement pass involves the application of this process with $R_i$ running once through all relators in the presentation.

Tietze processes function by working through these steps in some sensible order, guided by heuristics. Short eliminations and substring replacements reduce the total length of the presentation. Long eliminations can, and generally do, increase the length, often quite significantly. The substring searching component of substring replacements is by far the most time consuming part of Tietze processes, which is why we focus on it here.

## 2  Definition of the problem

For a presentation $G = \langle g_1, \ldots, g_d \mid R_1, \ldots, R_q \rangle$, we let $l_i$ denote the length of $R_i$ and $\Sigma$ denote the set of generators and their inverses $\{g_1, g_1^{-1}, \ldots, g_d, g_d^{-1}\}$ (the alphabet). **Rel** denotes the sequence of relators $\langle R_1, R_2, \ldots, R_{q-1}, R_q \rangle$, which is often kept sorted so that $l_1 \leq l_2 \leq \ldots \leq l_{q-1} \leq l_q$. We define the substring

searching (and replacement) problem that arises in Tietze transformation processes as:

Given **Rel** over $\Sigma$, for **any** two relators, $R_p$ and $R_t \in$ **Rel** with $l_p \leq l_t$, determine whether a common substring of length at least $\lceil (l_p + 1)/2 \rceil$ exists in equivalents of $R_p$ and $R_t$ and, if so, shorten $R_t$.

We use the following terminology. A **useful common substring** is a common substring of $R_i$ and $R_j$ with length greater than half the length of the shorter of $R_i$ and $R_j$. A common substring search between two **relators** $(R_p, R_t)$ and their equivalents is denoted by $ComStr(R_p, R_t)$, while the more usual common substring search between two **strings** $s_1$ and $s_2$ is denoted by $com\_substr(s_1, s_2)$. $R^\mathbf{i}$ denotes a string made from relator $R$ by rotating it $i$ positions right. The formal inverse of a string is obtained by reversing the string and inverting each symbol in the string (that is, replacing each $g_i$ by $g_i^{-1}$ and vice versa.) The equivalents of a relator $R$ which we consider are its rotations and their formal inverses. A **pass** is a substring replacement phase in which **each** pair of relators in **Rel** is considered once and only once for a $ComStr(R_p, R_t)$. If at least one of $R_p$ and $R_t$ has been changed since the previous $ComStr(R_p, R_t)$, then $ComStr(R_p, R_t)$ is a **necessary search** in the current pass, otherwise it is unnecessary since it is impossible that these relators have a useful common substring. We use $R_p$ to refer to a pattern relator and $R_t$ to a text relator.

We exemplify the performance of the various methods for substring searching applied to group presentations by considering some specific examples in detail. The performance gains demonstrated here typify the improvements achieved in this application area by these methods.

We study three applications, giving performance on presentations $\mathcal{J}$, $\mathcal{F}$ and $\mathcal{R}$. Presentation $\mathcal{J}$ is of the index 100 subgroup in the Janko simple group $J_2$, and comes from a subgroup presentation method (see [8]). It has 201 generators, 510 relators with longest relator of length 12, and total relator length 2,795. Presentation $\mathcal{F}$ is of the index 152 subgroup in the Fibonacci group $F(2, 9)$ and was obtained the same way. (It plays a crucial role in proving $F(2, 9)$ to be infinite, see [12, 15].) Presentation $\mathcal{F}$ has 153 generators, 304 relators with longest relator of length 13, and total relator length 2,119. Presentation $\mathcal{R}$ is for the restricted Burnside group $R(2, 5)$, a group of order $5^{34}$. It has 34 generators and 595 relators with longest relator of length 41, and total relator length 3,443. It was derived from a nilpotent quotient algorithm (see [10]).

## 3  Analysis of the problem

Algorithms and data structures for substring searching in various situations have been much studied, see [1] and [6, Chapter 7] for example. However the case considered here differs substantially from those covered there. Major distinguishing features of our situation are: all strings are (in effect) circular; formal inverses are (implicitly) present; many substrings are simultaneously sought; and the text is dynamic, changing very often. In this section we study features of our problem.

In our substring searching problem, each relator $R_i$ can be thought of as representing $2l_i$ strings: $l_i$ strings obtained by rotation; and another $l_i$ strings obtained by formal inversion. Thus **Rel**, which consists of $q$ relators, represents $2\sum_{i=1}^{q} l_i$ strings. The common substring search process for a pair of relators $R_p$ and $R_t$, $ComStr(R_p, R_t)$, can be concisely described in pseudocode in terms of $2l_p l_t$ $com\_substr$s (common substring searches for strings) as follows.

for $i := 0$ to $l_p - 1$ { for $j := 0$ to $l_t - 1$) {
$com\_substr(R_p^{\mathbf{i}}, R_t^{\mathbf{j}}); com\_substr((R_p^{\mathbf{i}})^{-1}, R_t^{\mathbf{j}}); $ } }

A pass consists of the $\binom{q}{2}$ choices of pairs of relators, that is, $q(q-1)/2$ $ComStr$s.

Since all rotations of $R_i$ are substrings of the string $R_i R_i$, it is not necessary to explicitly generate them all separately. A simple solution comes from relator extension. If we extend $R_p$ and $R_t$ by their initial $l_p/2$ symbols to obtain $R_p'$ and $R_t'$, then $ComStr(R_p, R_t) = \{com\_substr(R_p', R_t'); com\_substr(R_p^{-1'}, R_t'); \}$.

All of the relators which make up the presentation are used as patterns as well as texts. During the Tietze processing, they change frequently. Eliminations (short and long) and successful replacement passes make changes to relators. However, not all relators are changed between substring replacement passes. We use a two-level substring searching model, the skip level and the match level, to speed up the substring replacement passes.

At the skip level, unnecessary searches are identified and skipped. Early implementations of Tietze transformation programs compare each relator $R_i$ with every subsequent relator in the relator sequence in every pass. The idea here is that pairs of relators already searched are not searched again. Havas and Ollila [11] used change flags to avoid unnecessary searches. Here we improve on change flags by using a timestamp system. Two timestamp algorithms are given in the next section. Since many relators are not changed in a pass, this speeds up the whole process tremendously over the early methods, as the cost of timestamping is negligible.

The following practical results give the total number of relator pairs searched in equivalent Tietze processes on the given presentation. In the case of $\mathcal{J}$, a total of 6,693,105 searches were made with the early method. Change flags reduced this to 482,959 searches, further reduced to 351,253 by timestamps. (Only 2,376 of these were successful.) For $\mathcal{F}$, the corresponding figures are: 9,513,358 (early method); 832,689 (change flags); 585,383 (timestamps); and 2,739 (successful). Thus over 90% of the $ComStr$s which were done with the early method are skipped if we use change flags; timestamps provide a further 27% saving. Since the time used for handling change flags or timestamps is insignificant, searching time is similarly reduced.

At the match level, numerous variations are possible, with plenty of scope for improvement. This is because the successful search rate is very low, even when using timestamps, as illustrated above. The successful $ComStr$s (a useful common substring found) comprise only 0.68% and 0.47% of the necessary $ComStr$s for $\mathcal{J}$ and $\mathcal{F}$ respectively. Thus, if we can detect that there is no useful common substring quickly then a substantial time saving may be achieved.

## 4 Timestamps

In this section we present two timestamp algorithms and two theorems. One algorithm deals with "sorted relators" and the other with "unsorted relators". $R_i.T_p$ records the latest time when $R_i$ is used as a pattern. $R_i.T_s$ records the latest time when $R_i$ is changed in a $ComStr(R_p, R_i)$.

The following algorithm is applicable when **Rel** is kept sorted all the time (as in [9, 11]). Note that `Other operations` refers to elimination phases, which are also responsible for updating `R[*].Tp`, `R[*].Ts` and `NumRels`, as appropriate.

```
Initialize: R[*].Tp := -1; R[*].Ts := 0; timer := 1;
while (MoreSubstringSearchPass) {
  for p := 1 to NumRels-1 {
    for t := p+1 to NumRels {
      if (R[p].Tp <= R[t].Ts) {
        ComStr(R[p], R[t]);
        if R[t] changed {
          R[t].Tp := -1; R[t].Ts := timer; reorder R[t] in Rel;}
    } }
    R[p].Tp := timer; timer++;
  }
  Other operations;
  Compute MoreSubstringSearchPass;
}
```

**Theorem 1** *This algorithm performs all necessary searches and all searches the algorithm does are necessary.*

The proof is by detailed but straightforward analysis.

The following algorithm is used when **Rel** is not necessarily kept sorted within a pass (as in [18]). This is applicable if $R_k$ remains the $k$th relator in **Rel** throughout a pass, no matter whether it is changed or not. However at the beginning of each pass **Rel** is sorted.

```
Initialize: R[*].Tp := R[*].Ts := *;
while (MoreSubstringSearchPass) { TsLocal[*] := 0;
  for p := 1 to NumRels-1 {
    for t := p+1 to NumRels {
     if (R[t].len >= R[p].len and ( (TsLocal[p]+TsLocal[t])!=0
          or R[p].Tp > R[t].Tp or R[p].Tp <= R[t].Ts) )
       { ComStr(R[p],R[t]); if (R[t] changed) TsLocal[t] := p; }
    }
    R[p].Tp := p; R[p].Ts := TsLocal[p];
  }
  Sort(Rel);
  Other operations;
  Compute MoreSubstringSearchPass;
}
```

**Theorem 2** *This algorithm performs all necessary searches and all searches the algorithm does are necessary.*

Again the proof is by detailed but straightforward analysis.

## 5 Signatures, Fingerprints, Bloom filters, and Automata

In this section, we study methods for the match level. We can achieve efficiencies if we can detect unsuccessful searches early. There are two categories of string-matching algorithms: exact match algorithms such as brute-force, Knuth-Morris-Pratt, Boyer-Moore and Boyer-Moore derivatives, and automaton-based ones; and algorithms that initially allow errors, such as those of Harrison and Karp-Rabin. All of these are described in [6].

In spite of the theoretical worst case inferiority of brute force searching, its average case performance is linear in the length of the text being searched. Furthermore, Gonnet and Baeza-Yates [6, Table 7.4] show that it performs quite well in practice. In [9] a variant of brute force searching which enables a search for many strings simultaneously at no extra cost was used.

Thus, consider $R_p$ and $R_t$ with $l_p \leq l_t$. In order to shorten $R_t$ any useful common substring must have length greater than half the length of $R_p$. This means that it will contain either the first symbol of $R_p$ or a middle symbol, or the inverse of one of those. (Further, if $R_p$ is a nontrivial power, a useful substring must contain the first symbol or its inverse. Also, generators which are known from the presentation to be involutions are known to be their own inverses.) So the search starts by searching for one of at most four symbols as starting points in $R_t$. When such a match is found an attempt is made to extend the match circularly both backwards and forwards until it is long enough to be useful.

The first use of algorithms which allow errors to save time in this context was by Havas and Ollila [11], based on ideas of Harrison [7]. The speed up comes from the replacement of some time consuming substring searches by much faster tests which reveal that no useful match is possible. Strings are characterized by signatures. Fast signature generation and comparison often determines that one string cannot be a substring of another much more quickly than explicit string searching. Havas and Ollila used rotation and inversion invariant signatures well-suited to this context and present detailed performance results. This approach was reasonably successful, but signature computation and comparison is by no means free. Havas and Ollila concluded that change flags outperformed signatures, and this result extends to timestamps.

The Tietze procedures in GAP [18] use the Karp-Rabin algorithm in the substring searching part, combined with change flags. In this, strings are characterized by shorter entities called fingerprints. Efficiencies are achieved by manipulating fingerprints instead of the (possibly much longer) strings. The algorithm associates with each string $X$ a fingerprint $\phi(X)$. The search for a match initially compares short fingerprints. When a fingerprint match is found, an exact-match method (usually) has to be invoked to confirm whether the fingerprint match corresponds to an actual string match or is a false match. False matches may occur unless $\phi(X)$ is a one-to-one mapping, which would be unusual.

In GAP, a fingerprint (an integer) is associated with each minimal possibly-useful substring in each pattern relator and its equivalents. This means that $2l_p$ strings of length $\lceil(l_p+1)/2\rceil$ are characterized by $2l_p$ integers. Then fingerprints are computed for the $l_t$ length $\lceil(l_p+1)/2\rceil$ substrings of each text relator (and its rotations). In order to quickly search for fingerprint matches, the pattern fingerprints are stored in a type of hash table. The hash table is represented by a data structure called a Bloom filter [4], which reduces the amount of space required to contain the hash-coded information from that associated with conventional methods. The reduction in space is at the cost of some percentage of erroneous look-ups, which may be tolerable in some applications. The filter comprises a bit vector and several hash transformations.

The Bloom filter in GAP is organized so that fingerprints are represented by 3 bits, one bit in each of three bit-tables. Three hash functions compute table addresses for each fingerprint. When a match is found a brute force algorithm is then used to confirm whether it is an actual match, since both fingerprints and Bloom filters allow erroneous matches.

As long as the Bloom filter is reasonably loaded, these Tietze procedures work well. They are fast and space efficient. However presentation $\mathcal{R}$ causes problems. Almost all matches are false: 1,694,640 out of 1,716,314. Thus almost 99% of the matches are false, and there are only 21,674 actual matches. This leads to inordinate execution time, used in the brute force searches, and a total cpu time of about 10 hours on a fast Sparc machine.

Where do these false matches occur? Are they false fingerprint matches? Or is it in the Bloom filters?

We replaced the Bloom filters by an ordinary hash table (which is slower and uses more space). This revealed that only 283 out of 21,957 fingerprint matches are false, 1.2%. The total execution time is reduced to less than an hour, about 9% of that using Bloom filters (but the ordinary hash table uses 8 times more space). This indicates that almost all false matches occur in the Bloom filters. A further study reveals that the false matches mainly occur at a late stage of the processing, when the size of the alphabet is 4 and the length of pattern strings is over 500. (With "small" examples, such as $\mathcal{J}$ and $\mathcal{F}$ the time taken by the ordinary hash table is about twice that taken by 3-bit Bloom filters.) Using 4 bits instead of 3 bits to represent a fingerprint in Bloom filters (and using a similar hash function to produce the addresses of the fourth bit) reduces the number of false matches for $\mathcal{R}$ to 432,383 and the execution time by about a factor of three compared to the 3-bit filter. Again the false matches occur in the late stage when the size of the alphabet is 4, but even later, becoming frequent when the length of pattern strings is over 10,000 symbols.

Thus, fingerprints combined with Bloom filters provide an effective way of substring searching in this application. Except in the final stages of hard computations, when the filter may become overloaded, they are economical in both space and time. Alternative methods should be used in such final stages.

We have investigated the use of automaton-based string searching for this application. Automata have been successfully used to search for single pattern and multiple patterns [2, 3]. Perleberg [16] presented a longest substring (LS) searching algorithm based on automata. This algorithm requires another table *next length* in addition to *next state*. The *next length* table gives the maximum length of a substring that ends in the next state with the restriction that the next state follows the current state. Directly extending a single pattern match automaton to the LS problem would require $O(m^2|\Sigma|)$ space, $O(m^2|\Sigma| + m^3)$ preprocessing time, and $O(n)$ running time (where the pattern has length $m$ and the text length $n$). Perleberg's algorithm, by maintaining the *next length* table, only requires $O(m|\Sigma|)$ space, $O(m|\Sigma|+m^2)$ preprocessing time, and $O(n)$ running time.

Using relator extension as described in §3, we implemented Perleberg's algorithm at the match-level. For each pattern relator $R_p$, we build two automata, one for $R'_p$ and one for $R_p^{-1'}$. Even with a change to the heuristic strategy of Tietze processes to reduce the amount of substring searching, we found the automaton-based method to be slow. Thus, automaton searching takes 226 seconds for $\mathcal{J}$, compared with 37 seconds for a brute-force variant with change flags; for $\mathcal{F}$, it is 1,157 seconds as against 105. It takes too long to build the two automata for each pattern relator. For $\mathcal{J}$, there are 20,654 pattern relators in the whole run, for which automata construction takes 163 seconds, which is 72% of the total time. For $\mathcal{F}$, there are 26,431 pattern relators, and automata construction takes 75% of the the total time.

We can reduce the number of automata needed to one per search by using the equation $ComStr(R_p, R_t) = \{com\_substr(R'_p, R'_t);\ com\_substr(R'_p, R_t^{-1'});\ \}$. In this alternative, we replace a situation with two pattern strings and one text string by one with one pattern string and two text strings. Since we need one automaton per pattern, the time taken building automata is reduced by a factor of about two. Even though the preprocessing time is reduced, it is still far too much for our applications. The preprocessing time alone is still much more than the total time spent in the alternative method.

## 6  Conclusions

We have studied the substring searching component of presentation manipulation algorithms used in computational group theory. It differs from other string searching problems. We gave a formal definition of the problem and developed a two level searching model. We presented two timestamp algorithms at the first level and proved minimal-cover theorems associated with them. At the second level, we investigated methods based on signatures, fingerprints, Bloom filters, and automata. Detailed experiments revealed that different methods have advantages in different stages of the processes.